\documentclass[12pt,final]{amsart}
\usepackage{amsmath,amssymb}

\usepackage{cite}

\usepackage{color}
\usepackage{soul,graphicx}

\usepackage[color]{showkeys}
\definecolor{refkey}{rgb}{0,0,1}
\definecolor{labelkey}{rgb}{1,0,0}

\newcommand{\eq} [1] {\begin{equation}\label{#1}\quad}
\newcommand{\en} {\end{equation}}

\begin{document}

\begin{center}
		{\bf Stationary Processes, Wiener-Granger Causality, and Matrix Spectral Factorization}\\[5mm]
		Lasha Ephremidze

{\small			
 Kutaisi International University, Georgia

Razmadze Mathematical Institute of I. Javakhishvili Tbilisi State University
}
	\end{center}

\vskip+0.6cm

	{\small{\bf Abstract.} Granger causality has become an indispensable tool for analyzing causal relationships between time series. In this paper, we provide a detailed overview of its mathematical foundations, trace its historical development, and explore how recent computational advancements can enhance its application in various fields. We will not hesitate to present the proofs in full if they are simple and transparent. For more complex theorems on which we rely, we will provide supporting citations. We also discuss potential future directions for the method, particularly in the context of large-scale data analysis.
}
\vskip+0.2cm \noindent  {\small {\em  MSC:} 60G12, 47A68.

\vskip+0.2cm \noindent  {\small {\em  Keywords:} Stationary processes; Granger causality; matrix spectral factorization.

\section{Stationary processes in Hilbert spaces}

Let $H$ be a Hilbert space with a scalar product
$\langle\cdot,\cdot\rangle$ and a norm $\|\cdot\|$. Usually, we assume that $H$ is the space of square integrable functions $f$ on a probability space $(\Omega, \mathcal{F},P)$ with scalar product
$$
\langle f,g\rangle=\int\nolimits f\overline{g}\,dP.
$$
A sequence
$\{x_n\}_{n=-\infty}^\infty$ from $H$ is called stationary if
$$
\langle x_{n+m},x_n\rangle=c_m\;\;\text{ for each }n,m\in\mathbb{Z},
$$
i.e., the above scalar product is independent of $n$. The physical significance of stationarity lies in the invariance of the laws governing the process over time.

The sequence $c_m$ is positive definite, i.e.,
$$
\sum_{m=0}^N\sum_{n=0}^N c_{m-n}z_m\overline{z_n}\geq 0
$$
for any complex numbers $z_0,z_1,\ldots,z_N$. The latter inequality
holds since
$$
\sum\nolimits_{m=0}^N\sum\nolimits_{n=0}^N c_{m-n}z_m\overline{z_n}=
\left\|\sum\nolimits_{n=0}^N z_nx_n\right\|^2
$$
as one can check by direct computation. 

The positive definiteness of the sequence $c_m$, $m\in\mathbb{Z}$,
guarantees, by the classical theorem of Herglotz-Bochner \cite{FASP}, that there exists a positive measure on the unit circle $\mathbb{T}$, called the {\em spectral measure} of the process, such that 
\begin{equation}\label{2.1}
    c_m=\frac{1}{2\pi}\int\nolimits_{\mathbb{T}}t^{-m}\,d\mu(t)=
    \frac{1}{2\pi}\int\nolimits_0^{2\pi}e^{-im\theta}\,d\nu(\theta),
\end{equation}
where $\nu$ is the pullback measure of $\mu$ on $[0,2\pi)$
produced by the exponential map $\theta\mapsto e^{i\theta}$.

In the multidimensional case, we have a sequence
\begin{equation}\label{2.2}
    \{X_n\}_{n\in\mathbb{Z}}=\{(x_n^{(1)}, x_n^{(2)}, \ldots, x_n^{(d)})^T\}
\end{equation}
of $(d\times 1)$-columns of vectors from the Hilbert space, i.e,
$x_n^{(j)}\in H$, $n\in\mathbb{Z}$, $1\leq j\leq d$, and it is called a (jointly) stationary if
\begin{equation}\label{2.3}
   \langle x_{n+m}^{(k)}, x_n^{(l)}\rangle=c_m^{kl}.
\end{equation}
In high dimensions, we can use a matrix of (signed) measures
${\mathfrak{m}}=[\mu_{kl}]$ on $\mathbb{T}$ and it is called a positive definite
if the matrix
\begin{equation}\label{2.4}
   {\mathfrak{m}}(\Delta)=[\mu_{kl}(\Delta)]
\end{equation}
is positive definite for each Borel set $\Delta\in\mathbb{T}$, i.e.,
$$
(z_1,z_2,\ldots,z_d)\cdot\mathfrak{m}(\Delta)\cdot (z_1,z_2,\ldots,z_d)^T\geq
0$$
for each $\Delta\in\mathcal{B}$ and $ (z_1,z_2,\ldots,z_d)^T\in \mathbb{C}^{d\times 1}$. A
multidimensional extension to Bochner’s theorem assures that for a
jointly stationary sequence \eqref{2.2}, there exists a positive
definite matrix measure \eqref{2.4} such that \eqref{2.3} is equal to
$\frac{1}{2\pi}\int_\mathbb{T}t^{-m}\,\mu_{kl}(t)$ (see, e.g., \cite{Kolm41}). The latter equation
can be written in the matrix form
\begin{equation}\label{2.5}
 [X_{n+m}, X_n^T]=  \frac{1}{2\pi}\int\nolimits_{\mathbb{T}}t^{-m}\,d\mathfrak{m}(t)
\end{equation}
if we agree that $[X,Y]$, for
\begin{equation}\label{2.6}
   X=(x^{(1)}, x^{(2)}, \ldots, x^{(d)})^T \text{ and } Y=(y^{(1)}, y^{(2)}, \ldots, y^{(d)}) 
\end{equation}
is the matrix with $kl$th entry equal to $\langle x^{(k)},
y^{(l)}\rangle$ and the integration of a function with respect to a
matrix measure happens component-wise.

Note that, using the above notation, if $X$ is defined by \eqref{2.6}
and $(\alpha_1,\alpha_2,\ldots,\alpha_d)\in\mathbb{C}^{1\times d}$ is a
row vector of scalars, then
$$
\left\|\sum_{k=1}^d\alpha_kx^{(k)}\right\|^2=\left\langle
\sum_{k=1}^d\alpha_kx^{(k)},
\sum_{k=1}^d\alpha_kx^{(k)} \right\rangle=
(\alpha_1, \alpha_2,\ldots,\alpha_d)[X,X^T]
({\alpha_1},{\alpha_2},\ldots,{\alpha_d})^*.
$$
Furthermore, if $A_n=(\alpha_n^{(1)}, \alpha_n^{(2)},\ldots,\alpha_n^{(d)})\in \mathbb{C}^{1\times d}$, $n=0,1,\ldots,n$, and $\{X_n\}_{n\in\mathbb{Z}}$ is a stationary process \eqref{2.2}, then
\begin{gather}
\left\|\sum_{n=0}^N 
\sum_{k=1}^d\alpha_n^{(k)}x_n^{(k)}\right\|^2\!\!=\!
    \left\|\sum_{n=0}^N A_n X_n\right\|^2\!=\!\left\langle
\sum_{n=0}^N A_n X_n,
\sum_{n=0}^N A_n X_n \right\rangle\!=\!
\sum_{n=0}^N \sum_{m=0}^N  A_n[X_n,X_m^T]A_m^* \notag\\=
\sum_{n=0}^N \sum_{m=0}^N  A_n  \frac{1}{2\pi}\int\nolimits_{\mathbb{T}}t^{n-m}\,d\mathfrak{m}(t) A_m^*=
\frac{1}{2\pi}\int\nolimits_{\mathbb{T}} \left(\sum_{n=0}^NA_nt^n\right)\,d\mathfrak{m}(t) \left(\sum_{n=0}^NA_n^*t^{-n}\right). \label{sb7}
\end{gather}
The integration with respect to the matrix measure above should be understood in the usual measure-theoretic sense, with the customary ordering of operations for matrices.

These formulas for computing the norms of linear combinations of elements from \( H \) will be useful in determining the prediction errors discussed in the next section.

\section{Linear prediction theory of stationary processes}

Let $H_n\subset H$ be the minimal closed subspace spanned by
$\{x_n,x_{n-1},\ldots\}$. Following a
linear prediction theory of stationary processes created by Wiener \cite{Wiener49}
and Kolmogorov \cite{Kolm41}, projection of $x_L$, $L>0$, on $H_0$ is called an
$L$-lag ahead prediction of the process, namely
$$
\hat{x}_L:=Pr(x_L|H_0),
$$
and the norm
\begin{equation}\label{4.2}
e_x^{\{L\}}:=\|x_L-\hat{x}_L\|
\end{equation}
is called the $L$-lag prediction error.
Since the process is stationary, the quantity \eqref{4.2} is equal to 
$\|x_{n+L}-Pr(x_{n+L}|H_n)\|$ for each $n\in\mathbb{Z}$.

In the multidimensional case \eqref{2.2}, $H_n$ is the minimal closed
subspace spanned by $\{x_n^{(1)},\ldots, x_n^{(d)}, x_{n-1}^{(1)},\ldots,
x_{n-1}^{(d)},\ldots\}$, and the $L$-lag ahead prediction of $X_L$ is
taken to be
$$
\hat{X}_L=(\hat{x}_L^{(1)}, \hat{x}_L^{(2)}, \ldots,
\hat{x}_L^{(d)})^T, \;\text{ where }
\hat{x}_L^{(j}= Pr({x}_L^{(j)}|H_0), \;j=1,2,\ldots,d.
$$
The corresponding prediction error \eqref{4.2} is defined to be
$$
e_X^{\{L\}}=\|X_L-\hat{X}_L\|_{H^{d\times 1}}
:=\left(\sum\nolimits_{j=1}^d
\|x_L^{(j)}-\hat{x}_L^{(j)}\|^2\right)^{1/2}
$$
(to simplify notation, in what follows, for
$X=(x^{(1)},\ldots,x^{(d)})^T$, where $x^{(j)}\in H$, we assume
$\|X\|^2_{H^{d\times 1}}=\sum_{j=1}^d\|x^{(j)}\|^2$).

\section{Wiener-Granger causality}

Determining whether two time-varying processes are correlated is of fundamental importance in various branches of applied sciences. In most cases, this conclusion is drawn by analyzing time series, such as $\xi_n$ and $\zeta_n$, obtained from measurements of these processes, with several well-known standard techniques available for doing so. 
 However, in many important situations, we are interested in the causal relationships between the processes. Building on the existing theory of predicting stationary processes, Wiener introduced a brilliant new idea for measuring such relationships between time series \cite{Wien56}. The idea was quite simple. Wiener assumed that $\xi_n=x_n(\omega)$ and $\zeta_n=y_n(\omega)$ are realizations of corresponding stationary processes $x_n$ and $y_n$, and he proposed that the process generating $\zeta_n$ causes the process generating $\xi_n$, viz. $\zeta_n\to \xi_n$, if the prediction error
\begin{equation}\label{5.1}
e_{xy}=\inf\nolimits_{N,\alpha_0,\ldots,\alpha_N,\beta_0,\ldots,\beta_N}
    \left\| x_1-\sum\nolimits_{n=0}^N\alpha_nx_{-n}-\sum\nolimits_{n=0}^N\beta_ny_{-n}\right\|
\end{equation}
is significantly smaller than the error 
$$ e_{x}=\inf\nolimits_{N,\alpha_0,\ldots,\alpha_N}
    \left\| x_1-\sum\nolimits_{n=0}^N\alpha_nx_{-n}\right\|.$$
In other words, the inclusion of the second process in the prediction of the first significantly reduces the error norm (denoted as $e_{xy}<<e_x$ in such cases).

Nevertheless, Wiener's efforts went largely unnoticed because he described formula \eqref{5.1} using his own algorithm for matrix spectral factorization—a mathematical tool, which will be presented in the following sections—that appeared formidably complex. Granger \cite{Granger}, however, managed to express the same idea more simply using vector autoregressions and demonstrated several practical examples in economics, for which he gained both widespread recognition and a Nobel Prize in 2003.

Later, this idea was expanded in several directions. In particular, for vector observations $(x_n^{(1)}, \ldots, x_n^{(p)}, y_n^{(1)}, \ldots, y_n^{(q)})$ that are divided into two groups $X_n$ and $Y_n$, one can investigate whether one group causes the other. For this purpose, two quantities should be compared
\begin{equation}\label{6.1}
    e_{X}=\inf\nolimits_{N,A_0,\ldots,A_N}
    \left\| X_1-\sum\nolimits_{n=0}^N A_nX_{-n}\right\|_{H^{p\times 1}} \,,
\end{equation}
where $A_n$ are $p\times p$ matrices, $A_n\in \mathbb{C}^{p\times p}$, and
\begin{equation}\label{6.2}
    e_{XY}=\inf\nolimits_{N,B_0,\ldots,B_N}
    \left\| X_1-\sum\nolimits_{n=0}^N B_{n}(X_{-n}^T,Y_{-n}^T)^T\right\|_{H^{p\times 1}} \,,
\end{equation}
where $B_n\in\mathbb{C}^{p\times(p+q)}$, and we say that $Y_n\to X_n$ if $e_{XY}<<e_X.$

It may happen that the causal relationship between two processes, represented by the time series $x_n$  and $y_n$ or vector time series $X_n=
(x_n^{(1)}, x_n^{(2)}, \ldots, x_n^{(p)})$ and $Y_n=
(y_n^{(1)}, y_n^{(2)}, \ldots, y_n^{(q)})$, is delayed in time.
Therefore, there exists an approach \cite{BS19} that compares the quantities
\begin{gather}\label{s11}
 e_{x}^{\{L\}}=\inf\nolimits_{N,\alpha_0,\ldots,\alpha_N}
    \left\| x_L-\sum\nolimits_{n=0}^N\alpha_nx_{-n}\right\|,
\\ \label{s12}
 e_{xy}^{\{L\}} =\inf\nolimits_{N,\alpha_0,\ldots,\alpha_N,\beta_0,\ldots,\beta_N}
    \left\| x_L-\sum\nolimits_{n=0}^N\alpha_nx_{-n}-\sum\nolimits_{n=0}^N\beta_ny_{-n}\right\|
\end{gather}
and declares that there exists $L$-lag causal relation between $x_n$ and $y_n$, denoted $y_n\overset{\{L\}}{\to} x_n$, if $e_{xy}^{\{L\}}<<e_x^{\{L\}}$.

Accordingly, this idea can be generalized for vector observations, by taking $X_L$ instead of $X_1$ in formulas \eqref{6.1} and \eqref{6.2}. 

In what follows, we return to Wiener's approach and examine the role of matrix spectral factorization in the estimation of Granger causality.

\section{Spectral factorization}

Let $L_p(\mathbb{T})$, $p>0$, be the space of $p$-integrable functions on the unit circle $\mathbb{T}$ in the complex plane,
$$
L_p(\mathbb{T}):=\{f:\mathbb{T}\to\mathbb{C}\,|\, \int\nolimits_\mathbb{T}|f(t)|^p\,dt<\infty\},
$$
and let $\mathbb{H}_p(\mathbb{D})$, $p>0$, be the Hardy space of analytic functions $f$ in the unit disc $\mathbb{D}=\{z\in\mathbb{C}\,|\,|z|<1\}$ for which 
$$
\|f\|_{\mathbb{H}_p}^p=\sup\nolimits_{r<1}\int\nolimits_0^{2\pi} |f(re^{i\theta})|^p\,d\theta <\infty
$$
(see, e.g., \cite{Garn}, \cite{Koosis}). The boundary values function of $f\in \mathbb{H}_p$, i.e.
$$
f(e^{i\theta})=\lim\nolimits_{r\to 1}f(re^{i\theta}),
$$
exists a.e., belongs to $L_p(\mathbb{T})$, and it uniquely determines the function itself. Therefore, functions from $\mathbb{H}_p$ can be identified with their boundary values and we can assume that $\mathbb{H}_p(\mathbb{D})\subset L_p(\mathbb{T})$. Especially, we emphasize a natural characterization of $\mathbb{H}_p$, $p\geq 1$, in terms of the Fourier coefficients of boundary values:
$$
\mathbb{H}_p=\{f\in L_p\,|\, c_n\{f\}=0 \text{ for } n<0\},
$$
where 
$$
c_n\{f\}=\frac{1}{2\pi}\int\nolimits_0^{2\pi} f(e^{i\theta}) e^{-in\theta}\,d\theta, \;\;\;n\in\mathbb{Z},
$$
is the $n$-th Fourier coefficient of $f\in L_p$, $p\geq 1$.
(Whenever $f\in L_1(\mathbb{T})^{d\times d}$ is a matrix function, $c_n\{f\}\in \mathbb{C}^{d\times d}$.)

We also need the concept of {\em outer} analytic functions, which is somewhat technical. However, we will directly use the necessary properties without delving into the details. An analytic function $f\in \mathbb{H}_p$ is called outer, we use the notation $f\in \mathbb{H}_p^O$, if
\begin{equation}\label{7.5}
    |f(0)|=\exp\left(\frac{1}{2\pi}\int\nolimits_0^{2\pi} \log |f(e^{i\theta})|\,d\theta\right).
\end{equation}
The right-hand side of \eqref{7.5} is a maximal possible value of $|f(0)|$ in the class of functions $f$ from $\mathbb{H}_p$ with given absolute values on the boundary. In applications in signal processing, such functions are sometimes called also {\em optimal} or {\em minimal phase}.

The property of outer functions that will be used in the sequel is the following fact due to Beurling \cite[Th. II, 7.1]{Garn}. For any $f\in \mathbb{H}_p^O$, the set $\{fP\,|\, P\in\mathcal{P}\}$, where $\mathcal{P}$ is the set of polynomials, is dense in $\mathbb{H}_2$, i.e., for any $g\in\mathbb{H}_2$, there exists a sequence of polynomials $P_n(t)=\sum_{k=0}^n c_{nk}t^k$ such that
\begin{equation}\label{8.15}
\|g-fP_n\|_{\mathbb{H}_2}\to 0 \text{ as } n\to \infty.
\end{equation}

In the scalar case, the spectral factorization theorem asserts that for each non-negative $f\in L_1(\mathbb{T})$ which satisfies the Paley-Wiener condition,
\begin{equation}\label{8.2}
\int\nolimits_\mathbb{T}\log |f(t)|\,dt >-\infty,
\end{equation}
there exists $f_+\in \mathbb{H}_2^O$ such that
\begin{equation}\label{8.3}
f=|f_+|^2=f_+\overline{f_+} \text{ a.e. on }\mathbb{T}.
\end{equation}
There exists an explicit formula for computation of $f_+$ (\cite{Garn}, \cite{Koosis})
$$
f_+(z)=\exp\left(\frac{1}{4\pi}\int\nolimits_0^{2\pi} 
\frac{e^{i\theta}+z}{e^{i\theta}-z}
\log |f(e^{i\theta})|\,d\theta\right).
$$
and, therefore, it is assumed that $f_+$ can be constructed numerically for a given $f$.

In the multidimensional case, we have a positive definite $d\times d$ matrix function $S(t)$ with integrable entries,
\begin{equation}\label{9.0}
S=[s_{kl}]_{k,l=1}^d\in L_1(\mathbb{T})^{d\times d}
\end{equation}
and the Paley-Wiener condition has the form
\begin{equation}\label{9.1}
\int\nolimits_\mathbb{T}\log |\det S(t)|\,dt >-\infty.
\end{equation}
Under the circumstance, \eqref{9.0} can be factorized as
\begin{equation}\label{9.2}
S(t)=S_+(t)S_+^*(t),
\end{equation}
where, $S_+\in(\mathbb{H}_2(\mathbb{T}))^{d\times d}$ and $S_+^*$ is its Hermitian conjugate. In addition $S_+$ is {\em outer}, $S_+\in(\mathbb{H}_2(\mathbb{T}))^{d\times d}_O$, which by definition means  that $\det S_+\in \mathbb{H}_{2/d}^O$ (see \cite{EL10}). The generalization of Beorling's theorem asserts, in this case, that the set $$\{S_+(P_1,P_2,\ldots,P_d)^T\,|\, \text{ where } P_j\in\mathcal{P}\}$$ is dense in $(\mathbb{H}_2(\mathbb{T}))^{d\times 1}$, i.e., for each $(f_1,f_2,\ldots, f_d)^T$, where $f_j\in\mathbb{H}_2(\mathbb{T})$, there exists a sequence of vector polynomials 
$(P_n^{(1)}, P_n^{(2)}, \ldots, P_n^{(d)})^T\in \mathcal{P}^{d\times 1}$, where $P_n^{(j)}(t)=\sum_{k=0}^{N_n} c_{nk}^{(j)}t^n$ such that 
\begin{equation}\label{9.3}
\|(f_1,f_2,\ldots, f_d)^T-S_+(P_n^{(1)}, P_n^{(2)}, \ldots, P_n^{(d)})^T\|_{\mathbb{H}_2^{d\times 1}}\to 0 \;\;\text{ as } \;\;n\to 0.
\end{equation}

Numerical computation of \( S_+(t) \) for a given matrix function \eqref{9.2} has been a challenging problem for decades. Since Wiener’s original efforts \cite{Wie57}, \cite{Wie58} many authors have attempted to create a reliable algorithm for such factorization (see survey papers \cite{Kuc}, \cite{SayKai} and references therein). Nevertheless, the problem has remained far from a satisfactory solution, especially for large-scale singular matrices. The recent developments in the Janashia-Lagvilava method of matrix spectral factorization \cite{IEEE}, \cite{IEEE-2018} provide promising solutions that can overcome practical difficulties \cite{Ed22} .

\section{Computation of prediction errors}

Spectral factorization plays important role in computation of prediction errors $e_x, e_{xy}$, $e_x^{\{L\}}$, $e_{xy}^{\{L\}}$ and also $e_X, e_{XY}, e_X^{\{L\}}, e_{XY}^{\{L\}}$ for given stationary processes as demonstrated in this section. Without lose of much generality, we assume that the spectral measure $\mu$ in \eqref{2.1} or the matrix spectral measure $\mathfrak{m}=[\mu_{kl}]$ in \eqref{2.5} are absolutely continuous, i.e., $d\mu(t)=f(t)\,dt$, where $f\in L_1(\mathbb{T})$ and $d\mu_{kl}(t)=s_{kl}(t)\,dt$, where $s_{kl}\in L_1(\mathbb{T})$, and they satisfy the necessary and sufficient conditions \eqref{8.2} and \eqref{9.1}, respectively, for the existence of factorizations
\eqref{8.3} and \eqref{9.2}. Usually, spectral measures constructed from the observations satisfy these conditions, and the (matrix) functions $f$ and $S=[s_{kl}]$ are called, respectively, the power spectral density function of the process $x_n$ and power spectral density matrix of $X_n$.

We start with the scalar case and compute $e_x^{\{L\}}$ assuming $L\geq 1$:
\begin{gather*}
(e_x^{\{L\}})^2\!=\!\inf_{N,\alpha_0,\ldots,\alpha_N} \|x_L-\sum_{n=0}^N\alpha_nx_{-n}\|^2\!=\!\inf_{N,\alpha_0,\ldots,\alpha_N}\left\langle x_L-\sum_{n=0}^N\alpha_nx_{-n}, \;x_L-\sum_{n=0}^N\alpha_nx_{-n}  \right\rangle 
\displaybreak[0]\\
=\inf_{N,\alpha_0,\ldots,\alpha_N}\frac{1}{2\pi}\int_\mathbb{T}
(t^{-L}-\sum\nolimits_{n=0}^N\alpha_n t^n) (t^{L}-\sum\nolimits_{n=0}^N\overline{\alpha_n} t^{-n})
\,d\mu(t)
\displaybreak[0]\\=
\inf_{N,\alpha_0,\ldots,\alpha_N}\frac{1}{2\pi}\int_\mathbb{T}
(t^{-L}-\sum\nolimits_{n=0}^N\alpha_n t^n) f(t) (t^{L}-\sum\nolimits_{n=0}^N\overline{\alpha_n} t^{-n})\,dt \displaybreak[0]
\displaybreak[0]\\
=\inf_{N,\alpha_0,\ldots,\alpha_N}\frac{1}{2\pi}\int_\mathbb{T}
(t^{-L}-\sum\nolimits_{n=0}^N\alpha_n t^n) f_+(t)\overline{f_+(t)} (t^{L}-\sum\nolimits_{n=0}^N\overline{\alpha_n} t^{-n})\,dt
\displaybreak[0]\\=
\inf_{N,\alpha_0,\ldots,\alpha_N} \| (t^{-L}-\sum\nolimits_{n=0}^N\alpha_n t^n) f_+(t)\|_{L_2(\mathbb{T})}^2
\displaybreak[0]\\=
\inf_{N,\alpha_0,\ldots,\alpha_N} \left\| \sum\nolimits_{n=0}^{L-1} 
c_n\{f_+\}t^{n-L}+\sum\nolimits_{n=L}^{\infty} c_n\{f_+\}t^{n-L}-f_+(t)
\sum\nolimits_{n=0}^N\alpha_n t^n\right\|_{L_2(\mathbb{T})}^2
\displaybreak[0]\\=
\left\| \sum\nolimits_{n=0}^{L-1} c_n\{f_+\}t^{n-L} \right\|_{L_2(\mathbb{T})}^2\!+\!
\inf_{N,\alpha_0,\ldots,\alpha_N} \left\|\sum\nolimits_{n=L}^{\infty} c_n\{f_+\}t^{n-L}-f_+(t)
\sum\nolimits_{n=0}^N\alpha_n t^n\right\|_{L_2(\mathbb{T})}^2
\displaybreak[0]\\=
\left\| \sum\nolimits_{n=0}^{L-1} c_n\{f_+\}t^{n-L} \right\|_{L_2(\mathbb{T})}^2= \sum\nolimits_{n=0}^{L-1} |c_n\{f_+\}|^2.
\end{gather*}
The latter infimum is equal to $0$ because of \eqref{8.15}. Hence
\begin{equation}\label{10.10}
e_x^{\{L\}}= \left(\sum\nolimits_{n=0}^{L-1} |c_n\{f_+\}|^2\right)^{1/2}.
\end{equation}

Next we compute $e_{XY}^{\{L\}}$, which demonstrates a technique how to compute the quantities \eqref{5.1}-\eqref{s12}. We assume that $S(t)$ is a $(p+q)\times(p+q)$ power spectral density matrix of a stationary process 
$\{(X_{n}^T,Y_{n}^T)^T\}_{n\in\mathbb{Z}}$, where 
$$
X_n=(x_n^{(1)},x_n^{(2)},\ldots,x_n^{(p)})^T \text{ and }
Y_n=(y_n^{(1)},y_n^{(2)},\ldots,y_n^{(q)})^T,
$$
and \eqref{9.2} is its spectral factorization.

We use the following notation below: Let $e_d^{(1)}, e_d^{(2)},\ldots. e_d^{(d)}$ be the standard bases in the row space $\mathbb{C}^{1\times d}$. For a matrix $S$, let $S^{(j)}$ be its $j$th row, and $S^{(kj)}$ be its $kj$th entry.

We have (see \eqref{sb7})
\begin{gather*}
    \left(e_{XY}^{\{L\}}\right)^2=
    \inf_{N\geq 1, B_0,\ldots,B_N\in\mathbb{C}^{p\times(p+q)}}\left\|X_L-\sum\nolimits_{n=0}^N B_n\left(X_{-n}^T,Y_{-n}^T\right)^T\right\|^2_{H^{p\times 1}} 
    \displaybreak[0]\\=
     \inf_{N\geq 1, B_0,\ldots,B_N\in\mathbb{C}^{p\times(p+q)}} \sum_{j=1}^p \left\|X_L^{(j)}-\sum\nolimits_{n=0}^N B_n^{(j)}\left(X_{-n}^T,Y_{-n}^T\right)^T\right\|^2
     \displaybreak[0]\\
     =\sum_{j=1}^p \inf \left\langle X_L^{(j)}-\sum\nolimits_{n=0}^N B_n^{(j)}\left(X_{-n}^T,Y_{-n}^T\right)^T,\,
     X_L^{(j)}-\sum\nolimits_{n=0}^N B_n^{(j)}\left(X_{-n}^T,Y_{-n}^T\right)^T   \right\rangle
     \displaybreak[0]\\
     = \sum_{j=1}^p \inf \frac{1}{2\pi} \int_\mathbb{T}
     \left( e_{p+q}^{(j)} t^{-L}-\sum\nolimits_{n=0}^N B_n^{(j)}t^n\right)   S(t)
     \left( e_{p+q}^{(j)} t^{-L}-\sum\nolimits_{n=0}^N B_n^{(j)}t^n\right)^*dt
     \displaybreak[0]\\
      = \sum_{j=1}^p \inf \frac{1}{2\pi} \int_\mathbb{T}
     \left( e_{p+q}^{(j)} t^{-L}-\sum_{n=0}^N B_n^{(j)}t^n\right)
     S_+(t)S_+^*(t)
     \left( e_{p+q}^{(j)} t^{-L}-\sum_{n=0}^N B_n^{(j)}t^n\right)^*dt\\
     = \sum_{j=1}^p \inf \frac{1}{2\pi} \int_\mathbb{T}
     \left(\sum_{n=0}^{L-1}c_n\{S_+^{(j)}\}t^{n-L}+
     \sum_{n=L}^{\infty}c_n\{S_+^{(j)}\}t^{n-L}-
     \sum_{n=0}^N B_n^{(j)}t^n S_+(t)\right)
     \displaybreak[0]\\ 
     \times
      \left(\sum_{n=0}^{L-1}c_n\{S_+^{(j)}\}t^{n-L}+
     \sum_{n=L}^{\infty}c_n\{S_+^{(j)}\}t^{n-L}-
     \sum_{n=0}^N B_n^{(j)}t^n S_+(t)\right)^*dt
     \displaybreak[0]\\
     =\sum_{j=1}^p \frac{1}{2\pi} \int_\mathbb{T}
      \left(\sum\nolimits_{n=0}^{L-1}c_n\{S_+^{(j)}\}t^{n-L}\right)
      \left(\sum\nolimits_{n=0}^{L-1}c_n\{S_+^{(j)}\}t^{n-L}\right)^*dt
      \displaybreak[0]\\
      \left( \text{since } \inf\left\|\sum_{n=L}^{\infty}c_n\{S_+^{(j)}\}t^{n-L}-
     \sum_{n=0}^N B_n^{(j)}t^n S_+(t)\right\|_{\mathbb{H}_2^{1\times(p+q)}}=0 \text{ because of \eqref{9.3}} \right)\\
     \displaybreak[0]=\sum_{j=1}^p \frac{1}{2\pi} \int_\mathbb{T}
     \sum_{k=1}^{p+q} \left|\sum\nolimits_{n=0}^{L-1}c_n\{S_+^{(j)}\}t^{n-L}\right|^2 dt=
     \sum_{j=1}^p 
     \sum_{k=1}^{p+q} \sum_{n=0}^{L-1} \left|c_n\{S_+^{(jk)}\}\right|^2
          \end{gather*}
In a similar way, one can show that if $S(t)$ is a $d\times d$ power spectral density matrix of a stationary process \eqref{2.2} with spectral factor $S_+(t)$, then
\begin{equation}\label{s22}
    \left(e_{X}^{\{L\}}\right)^2=
 \sum_{j=1}^d 
     \sum_{k=1}^{d} \sum_{n=0}^{L-1} \left|c_n\{S_+^{(kj)}\}\right|^2.
\end{equation}

\section{Conclusions}
Although Hilbert spaces are often considered an abstract concept in mathematics, the theory described in this note can be directly applied in practice for causal analysis of multiple time series. Specifically, without explicitly referencing Hilbert spaces, the power spectral density matrix \eqref{9.0} can be constructed directly from time series observations, and as described above, this matrix encapsulates all the information about existing Wiener-Granger causal relations. There are several well-established methods for constructing power spectral densities \cite{Spec93}, and the method based on multitapers \cite{Thom82} has recently gained popularity among neuroscientists \cite{Dhamala} due to its efficiency. 

After constructing the power spectral density matrix, a novel algorithm for matrix spectral factorization \cite{IEEE-2018}, \cite{Ed22} can be employed to compute the corresponding coefficients and apply the relevant formulas \eqref{5.1}-\eqref{s12}, \eqref{10.10}-\eqref{s22} in practice.

\section{Acknowledgment}
The author was supported by the Shota Rustaveli National Science Foundation of Georgia (Project No. FR-18-18000).

\def\cprime{$'$}
\providecommand{\bysame}{\leavevmode\hbox to3em{\hrulefill}\thinspace}
\providecommand{\MR}{\relax\ifhmode\unskip\space\fi MR }
\providecommand{\MRhref}[2]{%
  \href{http://www.ams.org/mathscinet-getitem?mr=#1}{#2}
}
\providecommand{\href}[2]{#2}

\end{document}